\documentclass[conference]{IEEEtran}
\IEEEoverridecommandlockouts
\usepackage{cite}
\usepackage{amsmath,amssymb,amsfonts}
\usepackage{algorithmic}
\usepackage{graphicx}
\usepackage{textcomp}
\usepackage{soul,color}
\usepackage{xcolor}

\renewcommand{\b}{\boldsymbol}

\newcommand{\x}{\b{x}}
\newcommand{\lap}{\nabla^2}

\renewcommand{\L}{\mathcal{L}}
\newcommand{\T}{\mathsf{T}}
\newcommand{\R}{\mathbb{R}}

\def\BibTeX{{\rm B\kern-.05em{\sc i\kern-.025em b}\kern-.08em
    T\kern-.1667em\lower.7ex\hbox{E}\kern-.125emX}}
\begin{document}

\title{$p$-refined RBF-FD solution of a Poisson problem\\
{\footnotesize
\thanks{
    The authors would like to acknowledge the financial support of the ARRS
    research core funding No.~P2-0095 and the The World Federation
    of Scientists program.}
}}

\author{
    \IEEEauthorblockN{Mitja Jančič}
    \IEEEauthorblockA{
        ``Jožef Stefan'' International Postgraduate School \\
        Jamova cesta 39, 1000 Ljubljana, Slovenia \\
        and \\
        ``Jožef Stefan'' Institute \\
        Parallel and Distributed Systems Laboratory, \\
        Jamova cesta 39, 1000 Ljubljana, Slovenia \\
        Email: mitja.jancic@ijs.si}
    \and
    \IEEEauthorblockN{Jure Slak}
    \IEEEauthorblockA{``Jožef Stefan'' Institute \\
        Parallel and Distributed Systems Laboratory, \\
        Jamova cesta 39, 1000 Ljubljana, Slovenia \\
        Email: jure.slak@ijs.si}
    \and
    \IEEEauthorblockN{Gregor Kosec}
    \IEEEauthorblockA{
        \rule{11cm}{0px} \\[-13pt]
        ``Jožef Stefan'' Institute \\
        Parallel and Distributed Systems Laboratory, \\
        Jamova cesta 39, 1000 Ljubljana, Slovenia \\
        Email: gregor.kosec@ijs.si}
}

\maketitle

\begin{abstract}
    Local meshless methods obtain higher convergence rates when RBF approximations are augmented
    with monomials up to a given order. If the order of the approximation method is spatially
    variable, the numerical solution is said to be $p$-refined. In this work, we employ RBF-FD
    approximation method with polyharmonic splines augmented with monomials and study the numerical
    properties of $p$-refined solutions, such as convergence orders and execution time.
    To fully exploit the refinement advantages, the numerical performance is studied on a Poisson
    problem with a strong source within the domain.
\end{abstract}

\begin{IEEEkeywords}
    Meshless methods, $p$-refinement, RBF-FD, high order method
\end{IEEEkeywords}

\section{Introduction}
Meshless methods are becoming a strong alternative to mesh based methods, when numerical treatment
of partial differential equations is required. A strong advantage of meshless methods is that they
can operate on scattered nodes, contrary to mesh-based methods, that need a computationally
expensive mesh to operate. Many different meshless methods have been proposed so far, e.g.\
meshless Element Free Galerkin~\cite{belytschko1994element}, the Local
Petrov-Galerkin~\cite{atluri1998new}, h-p cloud method~\cite{duarte1996h} and others. In this paper
we use a method that generalizes the traditional
Finite Difference Method, called radial-basis-function-generated finite differences (RBF-FD). From
a historical point of view, RBF-FD was first mentioned by Tolstykh~\cite{tolstykh2003using}
in 2003 and has since been successfully used in a vast range of problems, e.g.\
convection-diffusion problems~\cite{chandhini2007local}, fluid flow problems~\cite{kosec2018local},
contact problems~\cite{slak2019adaptive}, scattering~\cite{slak2019high}, dynamic thermal rating of
power lines~\cite{maksic2019cooling}, etc.

The RBF-FD use RBFs to approximate the linear differential operators~\cite{bayona2010rbf}. Most of
the RBFs, like Hardy's multiquadrics or Gaussians, include shape
parameter~\cite{wendland2004scattered} that directly affects the stability of the approximation and
accuracy of the solution~\cite{schaback1995error}. To avoid shape parameter problems altogether,
Polyharmonic splines (PHS) have been proposed~\cite{bayona2017role}, however, PHS alone do not
guarantee convergent behavior. Therefore, RBFs are augmented with monomials up to given
order~\cite{bayona2017role}. The RBF part of the approximation takes care of the potential
ill-conditioning~\cite{flyer2016role}, while the polynomial part not only ensures convergent
behavior but also allows direct control over the convergence rate.

It has already been proven that having the control over the convergence rate is beneficial, when a
compromise between the accuracy of the solution and computational time is
needed~\cite{janvcivc2021monomial}. However, in this paper, we exploit the ability to control the
order of the approximation method to employ spatial variability of the approximating method order.
Such solution procedure refinement is also known as $p$-refinement~\cite{barros2004error} and is a
well known refinement procedure in the scope of finite element methods~\cite{babuska1992p},
where it also forms the basis of the highly successful $hp$-adaptive methods. In this
paper, convergence rates and computational times of $p$-refined solutions are studied. It is shown
that spatially variable order of the approximation method can notably reduce the computational time
and improve the convergence rate at the same time.

The rest of the paper is organized as follows: In section~\ref{sec:num} the main steps of solution
procedure are described, in section~\ref{sec:num_example}, the numerical example used to test the
numerical performance of $p$-refinement is presented, in section~\ref{sec:results}, results are
presented and finally, in section~\ref{sec:conclusions} conclusions are given and future work is
proposed.

\section{Numerical approximation}
\label{sec:num}
The solution procedure can be roughly divided into three steps. Using a dedicated node positioning
algorithm the domain is first discretized. Afterwards, the differential operators are approximated
in each node, resulting in stencil weights. Finally, the system of PDEs is discretized and,
therefore, transformed to a system of linear equations. The system is solved and its solution
stands for a numerical solution $u_h$ of the considered system of PDEs.

\subsection{Domain discretization}

In the first applications of meshless methods, researchers used existing mesh generators and
simply discarded the connectivity information to obtain the nodes~\cite{liu2003mesh}. However, such
procedure is computationally wasteful and conceptually wrong. Additionally, it can also be
problematic, since some authors reported that it failed
to produce nodal distributions of sufficient quality~\cite{shankar2018robust}.

\begin{figure}
    \centering
    \includegraphics[width=\linewidth]{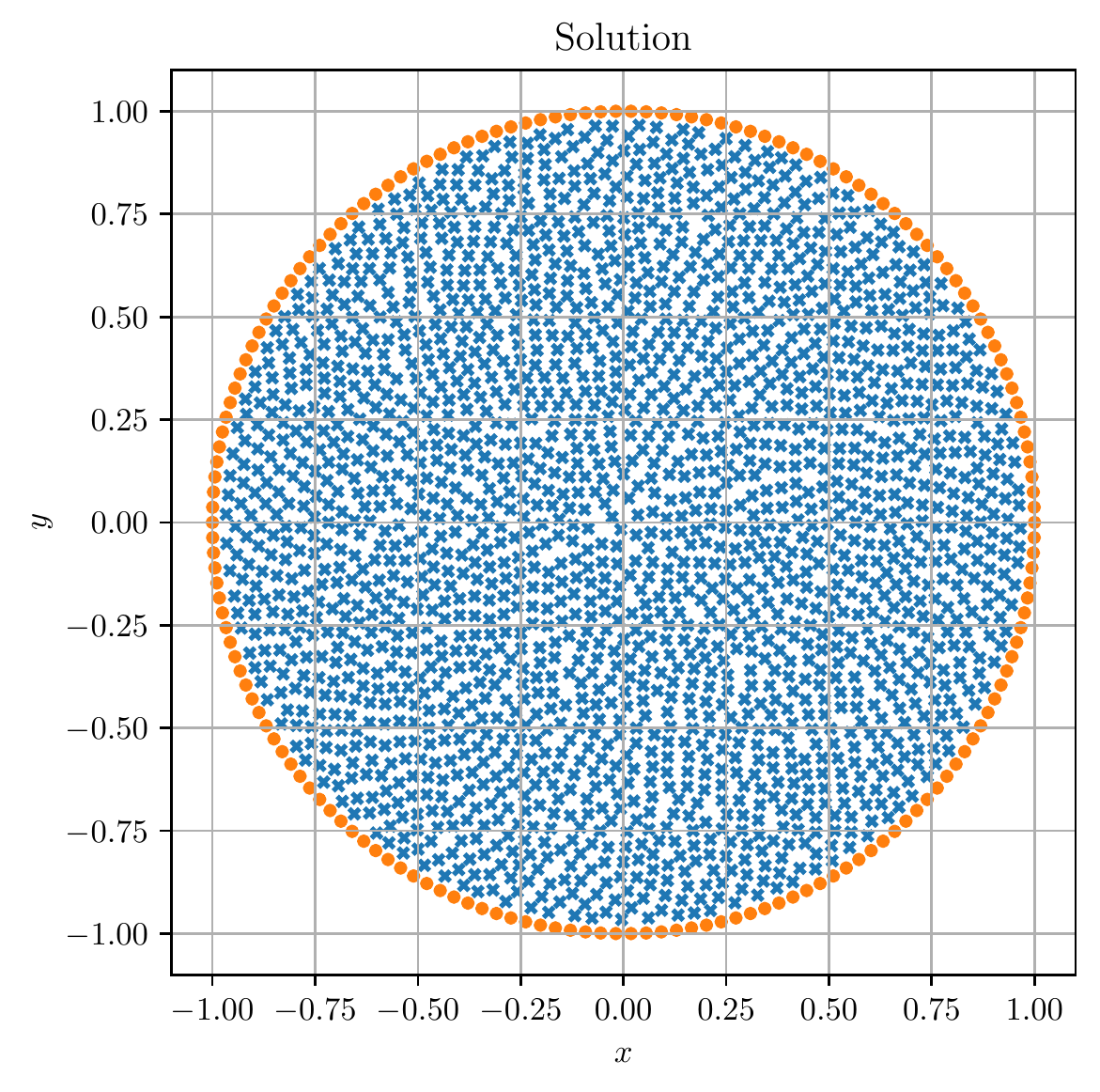}
    \caption{An example of circular domain populated with scattered nodes. Circles represent the nodes in the interior of the domain while the crosses are on its boundary.}
    \label{fig:domain_discretization}
\end{figure}

Researchers therefore soon started proposing dedicated node positioning algorithms. In this paper,
a dimension-independent node generation algorithm~\cite{slak2019generation} is used to
populate the domain with scattered nodes. The algorithm ensures a quasi-uniform internodal
distance $h$ as seen in figure~\ref{fig:domain_discretization}.

\subsection{Approximation of partial differential operators}
\label{sec:approx}

The behavior of a numerical method for solving systems of PDEs is defined by the approximation of
partial differential operators. In the scope of meshless methods, the approximation is done as
follows: Consider an operator $\L$ at a point $\x _c$. The approximation of $\L$ is sought
using the ansatz
\begin{equation}
    \label{eq:ansatz}
    (\L u)(\x _c) \approx \sum_{i = 1}^n w_i u(\x _i),
\end{equation}
where $w_i$ are the \emph{stencil weights}, $n$ is the \emph{stencil size} or \emph{support size},
and $u$ is the unknown function. To simplify the writing, the weights and function values are
assembled into vectors $\b w_\L(\x _c)$ and $\b u$ respectively. This notation allows us rewrite
the operator approximation~\eqref{eq:ansatz} in the form of a dot product between the two vectors
\begin{equation}
    \label{eq:ansatz_dot}
    (\L u)(\x _c) \approx \b w_\L(\x _c)^\T \b u.
\end{equation}

While the field values $\b u$ from equation~\eqref{eq:ansatz_dot} are considered as unknown, the
weights $\b w_\L(\x_c)$ need to be determined. To determine the weights, equality in
equation~\eqref{eq:ansatz_dot} is enforced for a given set of basis functions. In this paper we use
RBFs, denoted as $\phi _j$. The radially symmetric RBFS, centered at stencil nodes $\x _i$, can
be written in the form
\begin{equation}
    \phi _j (\x) = \phi(\left \| \x - \x _j \right \|),
\end{equation}
for a radial function $\phi$.
Each RBF corresponds to one linear equation
\begin{equation}
    \sum_{i = 1}^n \b w _i \phi _j(\x_i) = (\L \phi _j)(\x _c)
\end{equation}
with unknown weights $\b w_i$ and index $i$ running over all support nodes. Assembling these $n$
equations into matrix form, we obtain a system of linear equations
\begin{equation} \label{eq:rbf-system}
    \begin{bmatrix}
        \phi(\|\x_1 - \x_1\|) & \cdots & \phi(\|\x_n - \x_1\|) \\
        \vdots                & \ddots & \vdots                \\
        \phi(\|\x_1 - \x_n\|) & \cdots & \phi(\|\x_n - \x_n\|)
    \end{bmatrix}
    \begin{bmatrix}
        w_1 \\ \vdots \\ w_n
    \end{bmatrix}
    =
    \begin{bmatrix}
        l_{\phi, 1} \\
        \vdots      \\
        l_{\phi, n} \\
    \end{bmatrix}
\end{equation}
for
\begin{equation}
    l_{\phi, j} = (\L \phi(\|\x - \x_j\|))|_{\x = \x_c}.
\end{equation}
The system~\eqref{eq:rbf-system} is often compactly written as
\begin{equation}
    \b A\b w = \b l_\phi,
\end{equation}
where $\b A$ is symmetric and for some basis functions positive definite~\cite{buhmann2003radial}.

There are many different choices for RBFs. However, commonly used Hardy's multiquadrics or
Gaussians both depend on the shape parameter that directly affects the stability and accuracy of
the approximation~\cite{schaback1995error}. To avoid shape parameters entirely, Polyharmonic
splines (PHS) are used, defined as
\begin{equation}
    \phi(r) =
    \begin{cases}
        r^k,       & k \text{ odd}  \\
        r^k\log r, & k \text{ even}
    \end{cases},
\end{equation}
where $r$ denotes the Euclidean distance.

Note that using the PHS alone does not guarantee the convergence of local approximations from
equation~\eqref{eq:rbf-system}. Therefore, the approximation is additionally augmented with
monomials to omit the problems~\cite{bayona2017role}, which is done as follows.
Let $p_1, \dots, p_s$ be polynomials for which exactness of ansatz~\eqref{eq:ansatz_dot} is again
enforced. Monomials are often chosen up to a certain order $m$, resulting in $s=\binom{m+d}{m}$
monomials for a $d$-dimensional space.

The additional enforcement is introduced by extending the system~\eqref{eq:rbf-system} with the new
conditions
\begin{equation} \label{eq:rbf-system-aug}
    \begin{bmatrix}
        \b A    & \b P \\
        \b P^\T & \b 0
    \end{bmatrix}
    \begin{bmatrix}
        \b w \\ \b \lambda
    \end{bmatrix}
    =
    \begin{bmatrix}
        \b \ell_{\phi} \\ \b \ell_{p}
    \end{bmatrix}.
\end{equation}
Here
\begin{equation}
    \b P = \begin{bmatrix}
        p_1(\b x_1) & \cdots & p_s(\b x_1) \\
        \vdots      & \ddots & \vdots      \\
        p_1(\b x_n) & \cdots & p_s(\b x_n) \\
    \end{bmatrix}
\end{equation}
is a $n \times s$ matrix of polynomials evaluated at
stencil nodes and
\begin{equation}
    \b \ell_p = \begin{bmatrix}
        (\L p_1)|_{\b x=\b x _c} \\
        \vdots                     \\
        (\L p_s)|_{\b x=\b x_ c} \\
    \end{bmatrix}
\end{equation}
is the vector of values assembled by applying the considered operator $\L$ to
the polynomials at $\b x_c$.
Weights obtained by solving~\eqref{eq:rbf-system-aug}
are taken as approximations of $\L$ at $\b x_ c$ while additional unknowns $\b \lambda$, the Lagrange multipliers, are discarded.

The augmentation with monomials not only helps with convergence but also provides direct control
over the convergence rate, since the local approximation of the linear operator has the same order
as the basis used~\cite{bayona2017role}, while the RBF part handles the potential ill-conditioning
in purely polynomial approximations~\cite{flyer2016role}.

\subsection{PDE discretization}
Consider the boundary value problem with dirichlet boundary condition
\begin{align} \label{eq:bvp}
    \L u = f & \text{ in } \Omega,          \\
    u = g    & \text{ on } \partial \Omega.
\end{align}
The domain $\Omega$ is discretized with $N$ scattered nodes $\x_i$
with quasi-uniform internodal spacing $h$. Out of $N$ nodes, $N_i$ are in the interior and $N_d$ on the boundary $\partial \Omega$.

The stencils $\mathcal N(\b x_i)$ for each node $\b x_i$ are then selected. Commonly a single stencil constitutes of $n$ closest points, including the node itself. Choosing the right stencil size $n$ is far for trivial, however it has been recommended by Bayona~\cite{bayona2017role} to take at least $n=2\binom{m+d}{d}$ nodes.

In the next step, linear operator $\L$ is approximated at nodes $\x_i$, using the procedure
described in section~\ref{sec:approx}.  For each interior node $\x_i$, the equation $(\L u)(\x_i) =
f(\x_i)$ is approximated by a linear equation
\begin{equation} \label{eq:disc-eq}
    \b w_\L(\x_i)^\T \b u = \b f,
\end{equation}
where vectors $\b f$ and $\b u$ represent values of function $f$ and unknowns
$u$ in stencil nodes of $\x_i$. Similarly, for each Dirichlet boundary node $\x_i$, we obtain
the equation
\begin{equation} \label{eq:disc-dir}
    u_i = g(\x_i).
\end{equation}

All $N_i + N_d$ equations are assembled into a linear sparse system, with approximately $Nn$
nonzero elements. The solution $u_h$ of the system is a desired numerical approximation of $u$.
Note that using the spatially variable order of the approximation method can lead to a very
different number of nonzero elements in the linear sparse system.

\section{Numerical example}
\label{sec:num_example}
The behavior of the described solution procedure and its implementation is studied on a Poisson
problem with Dirichlet boundary condition. We aim to demonstrate and analyze the $p$-refined
solution procedure, where the order of the approximation method varies throughout the computational
nodes in the domain.

Governing equations are
\begin{align}
    \lap u (\x) & = f_{\text{lap}}(\x) & \text{ in } \Omega,         \\
    u (\x)      & = f(\x)       & \text{ on } \partial \Omega
\end{align}
where the domain $\Omega$ is a $d = 2$ dimensional circle
\begin{equation}
    \Omega = \left \{  \x \in \R^2, \left \| \x \right \| \leq \frac 3 2 \right \}.
\end{equation}
To fully exploit the advantages of $p$-refinement, the right hand side $f(\x)$ was chosen to have a relatively strong source within the domain at $\x_s = \b{\frac12}$, i.e.\
\begin{equation}
    f(\x) = \frac{1}{25 \left \| 4\x - \b 2 \right \|^ 2 + 1}.
\end{equation}
The Laplacian from $f$ can also be computed as
\begin{equation}
    f_{\text{lap}} (\x) = 3200\frac{25 \left \| 4\x - \b 2 \right \|^ 2}{f(\x)^3} -
    800\frac{d}{f(\x)^2} .
\end{equation}

The domain was filled with $N$ scattered nodes ranging from $N=1093$ to $N=978013$. The problem was
solved using RBF-FD with PHS augmented with monomials of order $m \in \left \{ 2, 4, 6 \right \}$.
Stencils for each node were selected by taking $n$ closest nodes where $n$ was determined as
recommended by Bayona~\cite{bayona2017role}
\begin{equation}
    n = 2\binom{m + d}{d}.
\end{equation}

\begin{figure}
    \centering
    \includegraphics[]{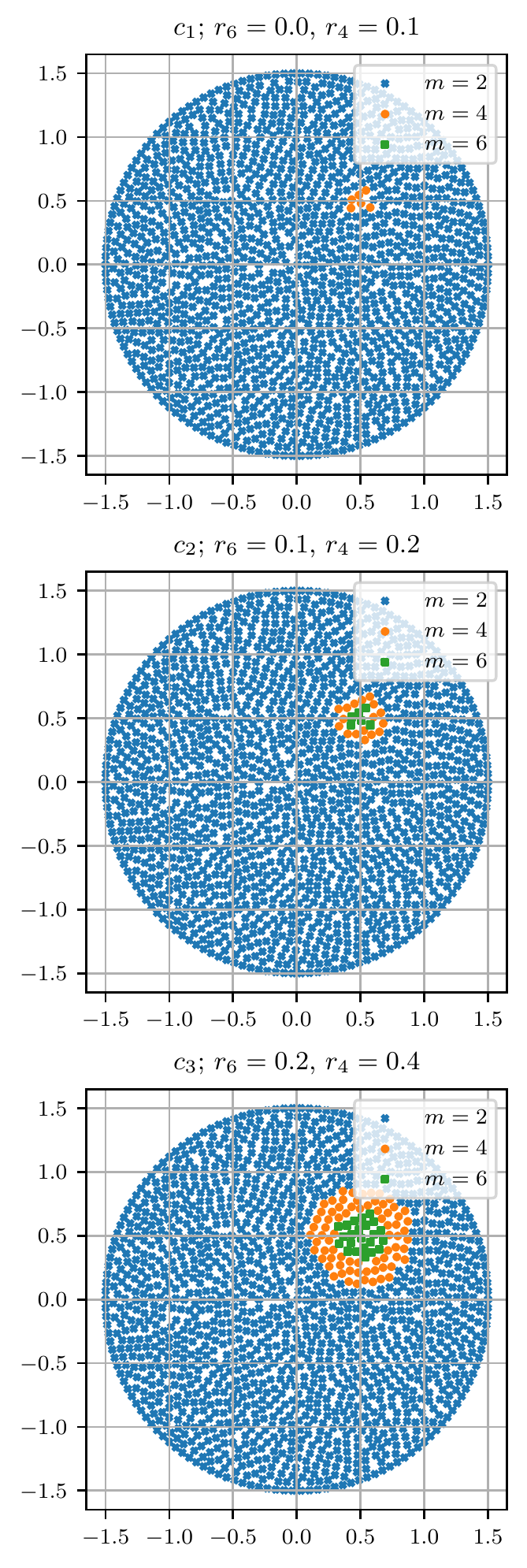}
    \caption{Three different stages of $p$-refinement used. Squares are used to mark the nodes where approximation is augmented with monomials of order $m=6$, circles for $m=4$ and crosses for $m=2$.}
    \label{fig:radiuses}
\end{figure}

To demonstrate the effect of $p$-refinement any combination of approximation orders $m$ can be
used. Naturally, to increase the overall convergence rate of the numerical solutions, the highest
approximation order is used where the numerical solution is expected to have the biggest error,
e.g.\ in the neighborhood of the strong source at $\x_s$. We define two radii $r_6$
and $r_4$ around source center $\x_s$. All computational nodes within the radius $r_6$, i.e.\
$\left \{ \x_i,\left \| \x_i - \x_s \right \| \leq r_6  \right \}$, are assigned with approximation
augmented with monomials of degree $m=6$, nodes within the annulus $\left \{ \x_i,r_6 < \left \|
\x_i - \x_s \right \| \leq r_4 \right \}$ are assigned approximation augmented with monomials of
degree $m=4$, while the remaining nodes are assigned approximation augmented with monomials of
degree $m=2$. So the order of the approximation method is spatially variable and can be compactly
written as
\begin{equation}
    m = \left\{\begin{matrix}
        6, & \left \| x_i - x_s \right \| \leq r_6     \\
        4, & r_6<\left \| x_i - x_s \right \| \leq r_4 \\
        2, & \text{otherwise.}
    \end{matrix}\right.
\end{equation}
Additionally, three different combinations $c_i$ of radii $r_6$ and $r_4$ have been used in this
paper
\begin{align}
    c_1 & = \left \{r_6= 0, r_4 = \frac{1}{10}\right \},           \\
    c_2 & = \left \{r_6= \frac{1}{10}, r_4 = \frac{1}{5}\right \} \text{ and} \\
    c_3 & = \left \{r_6= \frac{1}{5}, r_4 = \frac{2}{5}\right \}.
\end{align}
All three cases of spatially variable order of the approximation method are also shown in figure~\ref{fig:radiuses}.

\subsection{Error evaluation}
Closed form solution $u$ allows us to compute the accuracy of numerical solution $u_h$. In this paper, the error is evaluated in computational nodes with the infinity norm
\begin{equation}
    e_{\infty}  = \frac{\|u_h - u\|_\infty}{\|u\|_\infty}, \quad \|u\|_\infty = \max_{i=1, \ldots, N}
\end{equation}

The infinity norm is chosen as it is the strictest, but the authors also observed the same behavior
using 2-norm or 1-norm.

\section{Results}
\label{sec:results}
We compare the convergence rates of unrefined and $p$-refined numerical solution $u_h$ to the
problem from section~\ref{sec:num_example}.  Finally, we also study the effect of $p$-refinement on
computational times.

\subsection{Convergence rates}
\begin{figure}
    \centering
    \includegraphics[width=\linewidth]{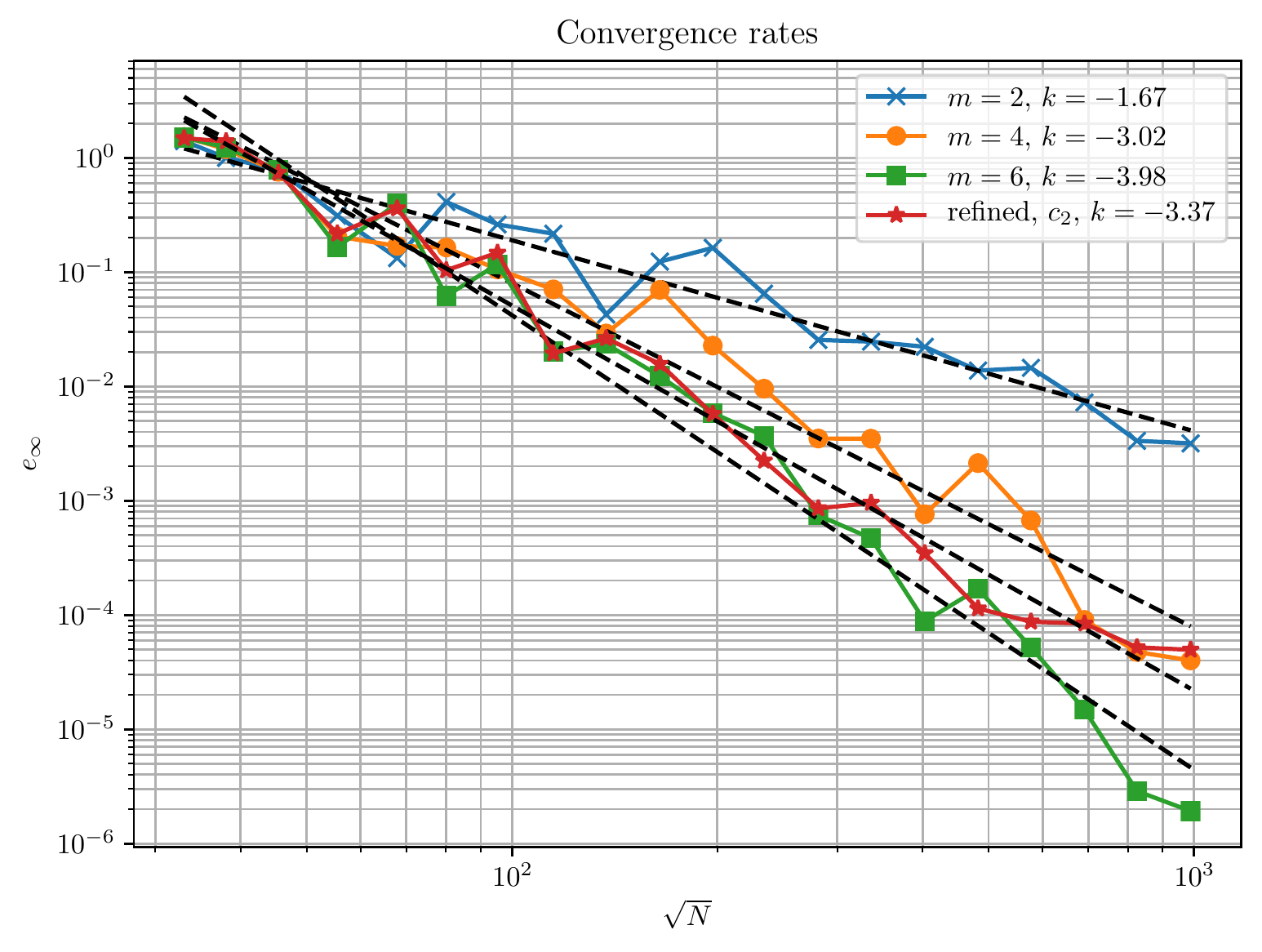}
    \caption{Convergence rates for different augmentation orders $m$ with respect to the number of nodes $N$.}
    \label{fig:conv}
\end{figure}

Convergence rates were estimated by computing the slope of the appropriate data subset. Selected
convergence rates are shown in figure~\ref{fig:conv}. We observe that the numerical
solutions converge for all chosen augmentation orders $m \in \left \{ 2, 4, 6 \right \}$. The
expected convergence rate of $O(h^m)$ is, however, not reached, but that is to be expected due to
the strong source within the computational domain.
The convergence curve of a $p$-refined solution for combination $c_2$ is
also added to the figure~\ref{fig:conv}. It is clear that the refined solution convergences at a
significantly better convergence rate compared to the convergence rate at $m=2$, regardless of the
fact that the majority of the stencils were still computed using monomials of order $m=2$. This
confirms our beliefs that the biggest contribution to the error comes from the strong source and
that the error can be, to some extent, mitigated by locally using a higher order method, i.e.\
$p$-refinement.

\begin{figure}
    \centering
    \includegraphics[width=\linewidth]{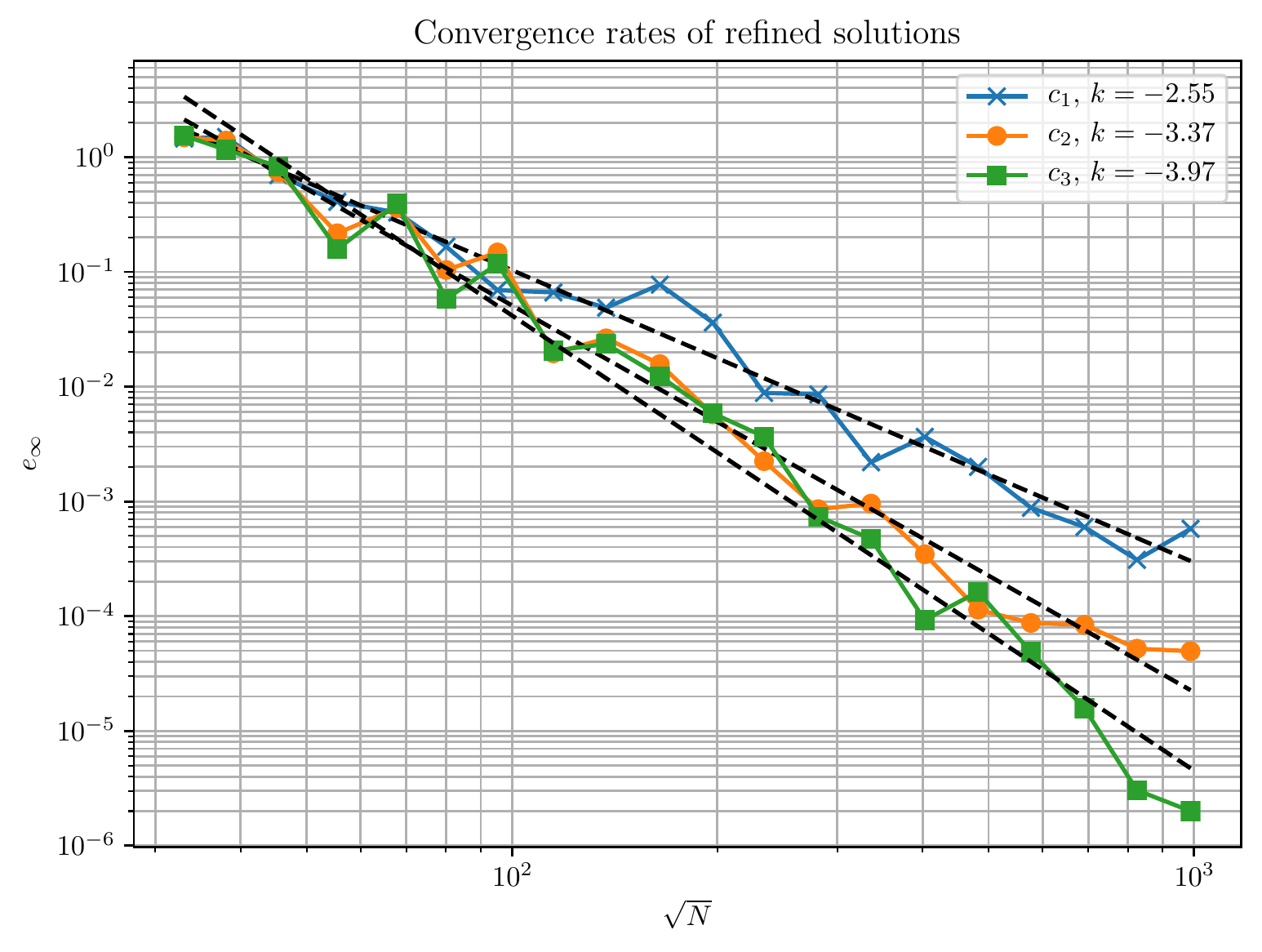}
    \caption{Refined convergence rates for different radius values combinations.}
    \label{fig:conv_ref}
\end{figure}

The effect of $p$-refinement is furthermore studied in figure~\ref{fig:conv_ref}, where convergence
rates of refined solutions for all combinations $\b c$ of radius values are shown. We observe how
the number of nodes used for higher order approximation affects the convergence rates. The
convergence rate for combination $c_3$ with the most higher order node stencils, is practically the
same as the convergence rate of unrefined solution with the highest augmented monomial $m=6$, even
though $m=6$ augmentation has only been applied to roughly 2 \% of all computational nodes and
$m=4$ to roughly 5 \%.

\subsection{Computation times}
In this section an overview of the total computational times is provided. All computations were
performed on a single core of a computer with
\texttt{AMD EPYC 7702 64-Core} processor and 512GB of DDR4 memory.
Code was compiled using \texttt{g++ (GCC) 9.3.0} for Linux with \texttt{-O3 -DNDEBUG} flags. The
sparse system is solved using the Pardiso solver on a single thread.

\begin{figure}
    \centering
    \includegraphics[width=\linewidth]{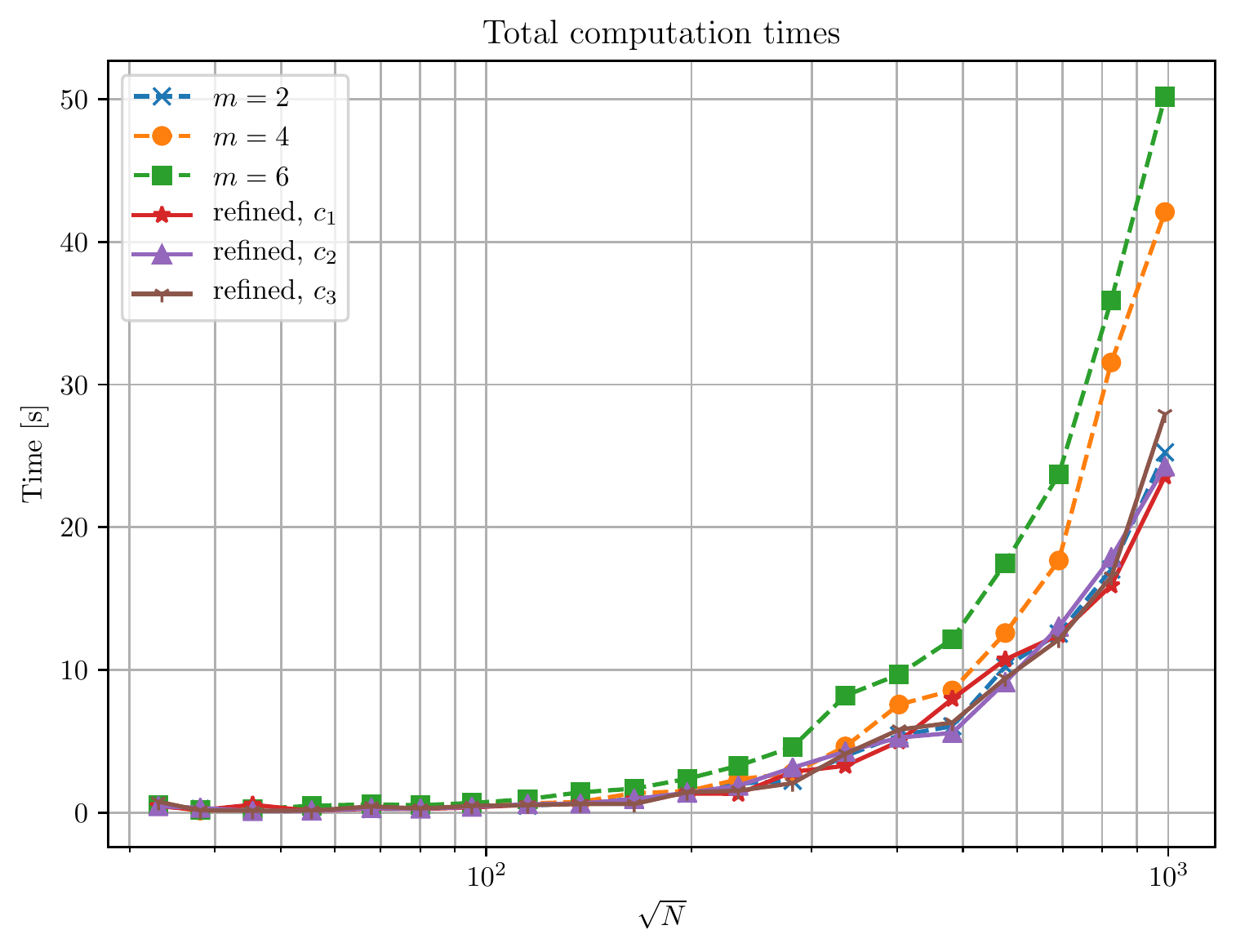}
    \caption{Total computational times with respect to the domain size $N$.}
    \label{fig:times}
\end{figure}

The total computational times are shown in figure~\ref{fig:times}, where the best run out of 5 is
selected. The total computational times of unrefined solutions (dashed lines) increases with the
monomial order $m$. This is expected, since the higher the order the larger the required stencil
size and consequently longer computational times. The computational times for all refined
solutions are similar to the unrefined solutions augmented with monomial order $m=2$, which
is also expected since the majority ($\approx$ 93 \%) of the nodes are assigned with augmentation
using monomials of degree $m=2$, however, results show that all refined solutions exhibit much
better convergence rates (see figures~\ref{fig:conv} and~\ref{fig:conv_ref}). This ultimately
means that employing the $p$-refinement enabled us to obtain significantly better convergence
behavior with little to no additional costs to execution time. Furthermore, refined solution
for $c_3$ combination with the largest radius values, measures the same convergence rate as
unrefined with $m=6$ ($k = -3.97$ vs.\ $k = -3.98$ respectively), but the former solution was
obtained approximately two times faster.

\section{Conclusions}
\label{sec:conclusions}
A $p$-refinement procedure the RBF-FD meshless method is presented, where the order of the local
approximation is spatially variable. We employed RBF-FD using PHS augmented with monomials of
different degrees to solve a Poisson problem with a strong source within the computational domain.
It is shown that $p$-refinement can improve the convergence rates at a very small cost to execution
time, and much faster, that using a method with a higher global order of convergence.

However, observations show that the $p$-refinement has its limitations. In some cases, specially
with local strong sources, the local description of the considered field is just not sufficient to
provide good local approximations of linear differential operators. Therefore a plan is to combine
the benefits of $p$-refinement with spatially variable nodal distributions, to provide better
approximations around critical areas within the domain. This is also known as $hp$-refinement, and
presents a major step towards $hp$-adaptivity.

\bibliographystyle{IEEEtran}
\bibliography{references}

\end{document}